\documentclass[11pt]{amsart}
\usepackage[T1]{fontenc}
\usepackage[latin1]{inputenc}
\usepackage{hyperref}
\usepackage{amscd}
\usepackage{epsfig}
\makeatletter
\theoremstyle{plain}
\newtheorem{thm}{Theorem}[section]
\numberwithin{equation}{section} %% Comment out for sequentially-numbered
\numberwithin{figure}{section} %% Comment out for sequentially-numbered
\theoremstyle{plain}    
\newtheorem{cor}[thm]{Corollary} %%Delete [thm] to re-start numbering
\newtheorem{lem}[thm]{Lemma} %%Delete [thm] to re-start numbering
\theoremstyle{plain}    
\newtheorem{prop}[thm]{Proposition} %%Delete [thm] to re-start numbering
\theoremstyle{remark}

\theoremstyle{remark}

\usepackage[arrow,matrix,curve,ps]{xy}
\xyoption{dvips}
\CompileMatrices

\renewcommand{\H}{\mathcal H}

\renewcommand{\O}{\mbox{$\mathcal{O}$}}
\renewcommand{\P}{\mbox{$\mathbb{P}$}}

\DeclareMathOperator{\birat}{bir}

\DeclareMathOperator{\Div}{Div}

\DeclareMathOperator{\Hom}{Hom}

\DeclareMathOperator{\locus}{locus}

\DeclareMathOperator{\RatCurves}{RatCurves}
\DeclareMathOperator{\red}{red}

\DeclareMathOperator{\Supp}{Supp}
\DeclareMathOperator{\Univ}{Univ}

\makeatother

\begin{document}

\title{Lines on contact manifolds II}

\date{\today }

\author{Stefan Kebekus}



\address{Stefan Kebekus, Institut für Mathematik, Universität
  Bayreuth, 95440 Bayreuth, Germany}

\email{stefan.kebekus@uni-bayreuth.de}

\thanks{The author gratefully acknowledges support by the
  Forschungsschwerpunkt ``Globale Methoden in der komplexen Analysis''
  of the Deutsche Forschungsgemeinschaft.}

\urladdr{http://btm8x5.mat.uni-bayreuth.de/$\sim$kebekus}

\maketitle
\tableofcontents

\section{Introduction}

Complex contact manifolds have recently received considerable
attention. Many of the newer publications, e.g.~\cite{Hwa97},
\cite{KPSW00}, \cite{Kebekus01} or \cite{Hon00}, approach contact
manifolds via the covering family of minimal rational curves, or via
geometric structures which are associated with the tangent vectors to
these curves. This short note furtheres the study of these curves.

We fix a family $H$ of minimal rational curves on a contact manifold
$X$. It is known that for any point $x\in X$, the subvariety
$\locus(H_x) \subset X$, which is covered by those curves which
contain $x$, is Legendrian. We will now study the deformations of
$\locus(H_x)$ which are generated by moving the base point.

As an application, we give a positive answer to a question of
J.M.~Hwang \cite{Hwa00} in the case of contact manifolds: for a
general point $x\in X$, the tangent map
$$
\begin{array}{rccc}
  \tau_x : & H_x & \to  & \P(T_X^*|_x)\\
  & \ell & \mapsto  & \P(T_\ell^*|_x)
\end{array}
$$
which maps a curve through $x$ to its tangent direction at $x$ is a
birational immersion. In other words, we show that a general choice of
a point $x\in X$ and a tangent direction $\vec v\in T_X|_x$ defines at
most a single minimal rational curve. The author believes that this is
a necessary step towards a full classification of contact manifolds.

We give a second application by showing that the normalization of
$\locus(H_x)$ is isomorphic to a projective cone.

\subsubsection*{Acknowledgement }

This paper was written up while the author enjoyed the hospitality of
the University of Washington at Seattle, the University of British
Columbia at Vancouver and Princeton University. The author would like
to thank K.~Behrend, J.~Kollár and S.~Kovács for the invitation.

\section{Setup}

Throughout the present paper, we maintain the assumptions and
notational conventions of the first part \cite{Kebekus01} of this
article. In particular, we refer to \cite{Kebekus01}, and the
references therein, for an introduction to contact manifolds and to
the parameter spaces which we will use freely thorughout.

For the reader's convenience, we will briefly recall the most
important conventions here. We assume that $X$ is a complex projective
manifold of dimension $\dim X = 2n+1$ which carries a contact
structure. This structure is given by a vector bundle sequence
$$
\begin{CD}
  0 @>>> F @>>> T_X @>{\theta}>> L @>>> 0
\end{CD}
$$
where $F$ is a subbundle of corank 1 and where the skew-symmetric
O'Neill-Tensor
$$
N : F\otimes F \to L
$$
which is associated with the Lie-Bracket is non-degenerate at every
point of $X$.

We will assume throughout that $X$ is not isomorphic to the projective
space $\P_{2n+1}$. It has been shown in \cite[sect.~2.3]{Kebekus01}
that this assumption implies that we can find a component $H \subset
\RatCurves^n(X)$ such that the intersection of $L$ with the curves
associated with $H$ is one. The space $H$ is therefore compact, and
for any point $x \in X$ we have a diagram
$$
\xymatrix{ {U_x} \ar[r]^(0.3){\iota_x}
  \ar[d]^{\pi_x}_{\txt{\scriptsize $\P_1$-bundle}}& {\locus(H_x) \subset X} \\
  {\tilde H_x} }
$$
where $\tilde H_x$ is the normalization the subfamily $H_x \subset
H$ of curves containing $x$ and $U_x$ is the pull-back of the
universal family $\Univ^{rc}(X)$. 

Throughout the paper we will constantly use the facts that that
$\tilde H_x$ is smooth (\cite[II.3.11.5]{K96}), that the subvariety
$\locus(H_x)$ which is covered by curves through $X$ is a Legendrian
subvariety of $X$ (\cite[prop.~4.1]{Kebekus01}) and that for a general
point $x \in X$ all curves $\ell \in H_x$ are smooth
(\cite[prop.~3.3]{Kebekus01}).

\section{Deformations of $\locus(H_x)$}

In order to study deformations of $\locus(H_x)$, it is useful to
consider sections $\sigma$ in the restriction of the tangent bundle
$T_X|_\ell$ to a minimal rational curve $\ell$. The following
proposition, which is a simple generalization of
\cite[prop.~3.1]{Kebekus01}, gives a convenient criterion which can be
used to show that $\sigma$ is contained in the restriction of $F$.

\begin{prop}\label{prop:sections}
  Let $x \in X$ be a general point, $f \in \Hom_{\birat}(\P_1, X,
  [0:1] \mapsto x)$ a morphism and let $\sigma \in H^0(\P_1,
  f^*(T_X))$ be a section such that $\sigma([0:1]) \in
  f^*(F)|_{[0:1]}$. 
  
  Then $\sigma$ is contained in $H^0(\P_1, f^*(F))$ if and only if
  $T_{\P_1}$ and $\sigma([0:1])$ are orthogonal with respect to the
  pull-back $f^*(N)$ of the O'Neill-tensor.
\end{prop}

\begin{proof}
  We know from \cite[thms.~II.3.11.5 and II.2.8]{K96} that the space
  $\Hom_{\birat}(\P_1,X)$ is smooth at $f$. Consequence: we can find
  an embedded unit disc $\Delta_\H \subset \Hom_{\birat}(\P_1,X)$,
  centered about $f$ such that $\sigma \in T_{\Delta_\H}|_f$ holds. In
  this situation we can apply \cite[prop.~3.1]{Kebekus01} to the
  family $\Delta_\H$, and the claim is shown.
\end{proof}

We will now employ this criterion in order to construct sections of
$L$ on the universal family over $H_x$.

\begin{lem}\label{lem:constructing_sections}
  If $x \in X$ is a general point, then there exists a natural vector
  space morphism
  $$
  s:T_X|_x \to H^0(U_x,\iota_x^*(L))
  $$
  such that the equality
  \begin{equation}\label{eq:prop_of_s}
    s(\vec v)|_{\sigma_\infty} = \iota_x^*(\theta(\vec v))    
  \end{equation}
  holds. 
\end{lem}

Note that the equation~(\ref{eq:prop_of_s}) makes sense because the
restricted line bundle $\iota_x^*(L)|_{\sigma_\infty}$ is trivial.

\begin{proof}
  If $\ell \in H_x$ is any line, then it follows immediately from
  proposition~\ref{prop:sections} that the kernel of the evaluation
  map
  $$
  e:H^0(\ell, T_X|_\ell) \to T_X|_x
  $$
  is contained in $H^0(\ell,F|_\ell)$ ---for this, note that the
  the vector $0\in T_X|_x$ is contained in $F$ and is perpendicular to
  any other vector. Recall from \cite[lem.~3.5]{Kebekus01} that the
  vector bundle $T_X|_\ell$ is globally generated. The evaluation map
  $e$ is therefore surjective, and we obtain a vector space morphism
  $$
  s_\ell : T_X|_x \to H^0(\ell, L|_\ell).
  $$
  Since $\tilde H_x$ is smooth, it is elementary to see that for
  any tangent vector $\vec v \in T_X|_x$, the union of the sections
  $(s_\ell(\vec v))_{\ell \in H_x}$ gives a section $s(\vec v)\in
  H^0(U_x,\iota_x^*(L))$.
\end{proof}

For those tangent vectors $\vec v$ which are contained in $F|_x$, we
are able to describe the section $s(\vec v)$ in greater detail. 

\begin{lem}\label{lem:special_section}
  If $\sigma_\infty \subset U_x$ is the section which is contracted by
  $\iota_x$, and if $\vec v\in F|_x \setminus \{ 0 \}$, then the
  support of the divisor $\Div(s(\vec v))$ which is associated with
  $s(\vec v)$ is given as
  $$
  \Supp(\Div(s(\vec v))) = \sigma_\infty \cup (\tau_x \circ
  \pi_x)^{-1}(D),
  $$
  where $D\in |\O_{\P(F|_x^*)}(1)|$ is the union of those tangent
  directions which are perpendicular to $\vec v$.
\end{lem}

\begin{proof}
  Since the pull-back $\iota_x^*(L)$ intersects $\pi_x$-fibers with
  multiplicity one, and since the section $s(\vec v)$ will always
  vanish on $\sigma_\infty$, it is clear that the support $\Supp(\Div(s(\vec v)))$ must
  be of the form 
  $$
  \Supp(\Div(s(\vec v))) = \sigma_\infty \cup  \pi_x^{-1} (D').
  $$
  The proof is finished if we show that any section $\sigma \in
  H^0(\ell, T_X|_\ell)$ with $\sigma(x) = \vec v$ is contained in
  $H^0(\ell,F|_\ell)$ if and only if the orthogonality holds. That,
  however, is exactly the statement of
  proposition~\ref{prop:sections}.
\end{proof}

\begin{cor}
  The morphism $s: T_X|_x \to H^0(U_x, \iota_x^*(L))$ is injective
  and the linear system $|\iota_x^*(L)|$ is basepoint-free.
\end{cor}

\begin{proof}
  By lemma~\ref{lem:special_section}, the restriction $s|_{F_x}$ is
  injective and the base locus of the linear sub-system $(s(\vec
  v))_{\vec v\in F|_x}$ is exactly the contracted section
  $\sigma_\infty$. The claim is therefore shown if we note that for
  every tangent vector $\vec v \in T_X|_x$, which is not contained in
  $F|_x$, the section $s(\vec v)$ does not vanish on
  $\sigma_\infty$. For this, we refer to equation~(\ref{eq:prop_of_s})
  of lemma~\ref{lem:constructing_sections} above.
\end{proof}

\section{Birationality of the tangent map}

We apply the results of the previous section to show that the tangent
map $\tau_x$ is generically injective. Together with the results of
\cite[thm.~3.4]{Keb00a} and \cite[cor.~3.6]{Kebekus01}, this implies
that $\tau_x$ is a finite birational immersion.

We start the proof by studying curves which intersect $\locus(H_x)$
tangentially. For the formulation of the proposition, recall the fact
\cite[thm.~3.5]{Keb00a} that for a general point $y\in \locus(H_x)$
there exists a unique curve in $\ell_{x,y} \in H_x$ which contains
both $x$ and $y$. The evaluation morphism $\iota_x$ is therefore
birational.

\begin{prop}\label{prop:intersect}
  If $y\in \locus(H_x)$ is a general point, and if $\ell_y \in H$ is a
  curve which intersects $\ell_{x,y}$ tangentially in $y$, then
  $\ell_y = \ell_{x,y}$.
\end{prop}

The proof will involve the infinitesimal description of the
$\Hom$-scheme and of the universal morphism $\mu :
\Hom_{\birat}(\P_1,X)\times \P_1 \to X$. The reader might want to
consider \cite[sect.~2.2]{Kebekus01} or \cite{Kebekus-Habil} for a
brief statement of the basic facts and the associated notation.

\begin{proof}
  \begin{figure}
    $$
    \xymatrix{
      { \ \begin{picture}(4,4)
          \put(0.0, 0.0){\epsfysize 3cm \epsffile{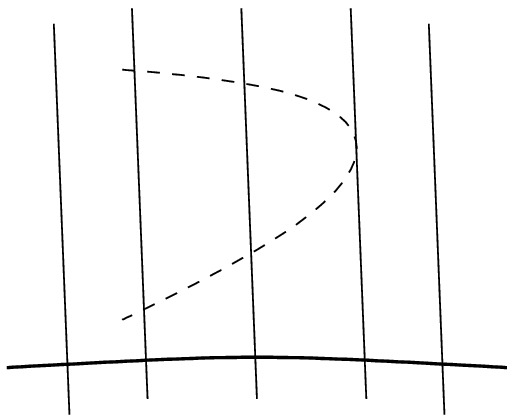}}
          \put(2.0, 3.0){\small $U_x$}
          \put(0.7, 0.9){\small $\tilde \ell_y$}
          \put(2.5, 1.9){\small $\bullet$}
          \put(2.4, -0.1){\small $\tilde \ell_{x,y}$}
          \put(-0.6, 0.5){\small $\iota_x^{-1}(x)$}
        \end{picture}
        } \ar[r]^{\iota_x} & { \begin{picture}(4,4)
          \put(0.0, 1.0){\rotatebox{-40}{\epsfysize 4cm \epsffile{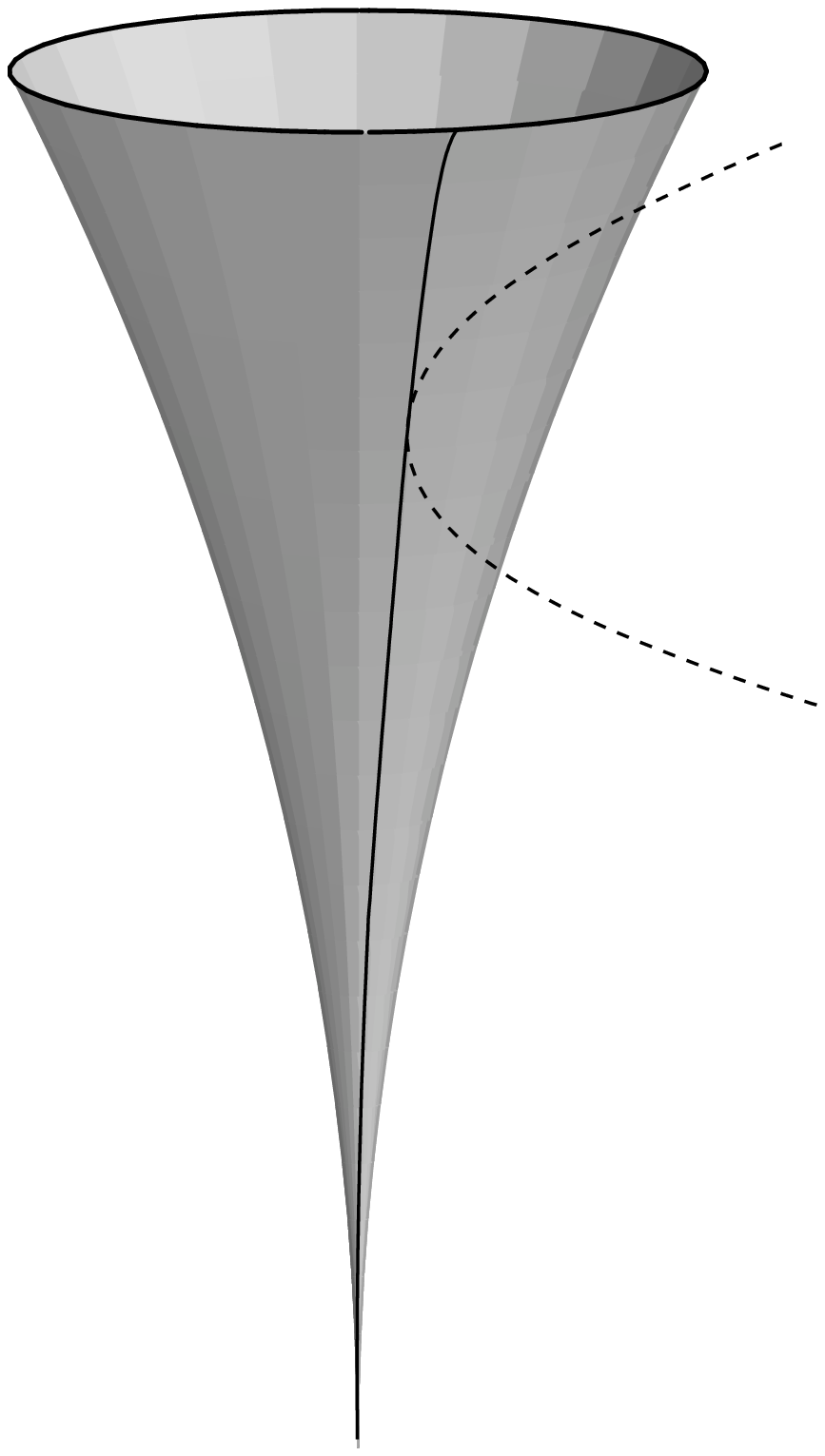}}}
          \put(4.0, 3.6){\small $\locus(H_x)$}
          \put(4.1, 2.2){\small $\ell_y$}
          \put(2.35, 2.25){\small $y \, \bullet$}
          \put(0.57, 0.36){\small $x \, \bullet$}
        \end{picture}
        }
      }
    $$
    \caption{Proof of proposition~\ref{prop:intersect}}
  \end{figure}
  As a first step, we claim that the curve $\ell_y$ is contained in
  $\locus(H_x)$.

  \subsubsection*{Proof of the claim} 
  
  As $y$ was chosen generally, we can find a birational morphism $f_y
  : \P_1 \to \ell_y$ with $f_y([0:1])=y$ such that the reduced scheme
  $\Hom_{\birat}(\P_1,X)_{\red}$ is smooth at $f_y$ and such that
  tangent map of the restricted universal morphism
  $\mu|_{\Hom_{\birat}(\P_1,X)_{\red}\times\{[0:1]\}}$ has maximal
  rank at $f$. This enables us to find an analytic open subset $U_x
  \subset \locus(H_x)$ and a lifting $U_\H \subset
  \Hom_{\birat}(\P_1,X)$ with curves which are tangent to $\ell_{z,x}$
  for points $z\in U_x$. More precisely, we can find an embedded
  polycylinder $U_\H$ such that
  \begin{enumerate}
  \item The restricted universal morphism
    $$
    \begin{array}{rccc}
      \mu|_{U_\H\times \{[0:1]\}} : & U_\H & \to  & X\\
      & f_z & \mapsto  & f_z([0:1])
    \end{array}
    $$
    induces a morphism $U_\H \to U$.
    
  \item For every morphism $f_z \in U_\H$, the curve $\ell_z :=
    f_z(\P_1)$ meets $\ell_{x,z}$ tangentially in $z$.
  \end{enumerate}
  We will now apply proposition~\ref{prop:sections} to sections
  $\sigma \in H^0(\ell_z,T_X|_{\ell_z})$ which come from the
  deformation family $U_\H$. It is clear from the construction that if
  $$
  \sigma \in T_{U_\H}|_{f_z} \subset H^0(\ell_z,T_X|_{\ell_z})
  $$
  is any tangent vector, then
  $$
  \sigma([0:1]) \in T_{\locus(H_x)}|_z.
  $$
  
  Since $\locus(H_x)$ is known to be Legendrian, the tangent space
  $T_{\locus(H_x)}|_z$ will be isotropic with respect to the
  O'Neill-tensor $N$; the restriction $N|_{T_{\locus(H_x)}|_z}$ is
  identically zero:
  $$
  T_{\ell_z}|_z \perp T_{\locus(H_x)}|_z.     
  $$
  In particular, the tangent vector $\sigma([0:1])$ is
  perpendicular to $T_{\ell_z}|_z$. Applying
  proposition~\ref{prop:sections}, we see that the tangent space
  $T_{U_\H}|_{f_z}$ is associated with sections in
  $H^0(\ell_z,F|_{\ell_z})$. This implies that the image
  $\mu(U_\H\times \P_1)$ must be $F$-integral ---see
  \cite[prop.~II.3.4]{K96}. Since, on the other hand, $\mu(U_\H\times
  \P_1)$ already contains the Legendrian submanifold $U$, it turns out
  that $\mu(U_\H\times \P_1)$ is contained in $\locus(H_x)$ so that
  the claim $\ell_y \subset \locus(H_x)$ follows.
  
  \subsubsection*{Application of the claim} 
  As a next step we consider the strict transform $\tilde \ell_y$ of
  $\ell_y$ in the universal family $U_x$. We will now argue by
  contradiction, assume that $\tilde \ell_y$ is not a fiber of the map
  $\pi_x: U_x\to \tilde H_x$, and derive a contradiction. For this,
  note that the curve $\tilde \ell_y$ intersects the strict transform
  $\tilde \ell_{x,y}$ of $\ell_{x,y}$, which is a $\pi_x$-fiber,
  transversally in the preimage of $y$. We claim that this is
  impossible; a contradiction is thus reached.
  
  In order to see this contradiction, observe that $\iota_x^*(L)$
  intersects the curve $\tilde \ell_y$ with multiplicity one; this is
  because $y$ is general and the morphism $\iota_x:U_x\to \locus(H_x)$
  is birational. On the other hand, by lemma~\ref{lem:special_section}
  we can always find a vector $\vec v \in T_X|_x$ such that the
  associated section $s(\vec v) \in H^0(U_x,\iota_x^*(L))$ vanishes on
  the fiber $\tilde \ell_{x,y}$, but not on the curve $\tilde \ell_y$.
  Thus, the curve $\tilde \ell_y$ intersects the divisor $\Div(s(\vec
  v))$ tangentially, and the intersection number cannot be one. This
  ends the proof of proposition~\ref{prop:intersect}.
\end{proof}

\begin{cor}
  The tangential map $\tau_x$ is generically injective, that is,
  $\tau_x$ is birational onto its image.
\end{cor}

\begin{proof}
  Using the notation of above, if $\ell_{xy}\in H$ is the unique curve
  containing both $x$ and $y$, and if $\ell_y \in H$ is a curve which
  contains $y$, then proposition~\ref{prop:intersect} shows that
  $\ell_{x,y}$ and $\ell_y$ intersect transversally. Since $y$ was
  generically chosen, the claim follows.
\end{proof}

\section{The normalization of $\locus (H_x)$}

We will now show that the normalization of $\locus(H_x)$ is isomorphic
to a cone. 

\begin{prop}
  If $x$ is a general point, then the normalization
  $\widetilde{\locus(H_x)}$ of $\locus(H_x)$ is a cone.
\end{prop}

\begin{proof}
  Since all curves which are associated with points in $H_x$ are
  smooth, and since $\tilde H_x$ is smooth, the preimage
  $\sigma_\infty := \iota_x^{-1}(x) \subset U_x$ is a section over
  $\tilde H_x$ which can be contracted to a point. By definition,
  $\widetilde{\locus(H_x)}$ is a cone if we can show that
  \begin{enumerate}
  \item $\widetilde{\locus(H_x)}$ is isomorphic to the image of $U_x$
    under the contraction of $\sigma_\infty$.

  \item There exists another section $\sigma_0 \subset U_x$ which is
    disjoint from $\sigma_\infty$.
  \end{enumerate}
  See the introductory remarks in the paper \cite{Wah83} for more
  information on this.
  
  In order to show property (1), recall Mori's Bend-and-Break argument
  \cite[II.5]{K96} which asserts that $\iota_x|_{U_x\setminus
    \sigma_\infty}$ is finite. Likewise, recall from
  \cite[thm.~3.6]{Keb00a} that $\iota_x$ is birational. Consequence:
  the Stein factorization of $\iota_x$ yields a decomposition as
  follows:
  $$
  \xymatrix{ {U_x} \ar[rr]_(.425){\txt{\scriptsize contr.~of
        $\sigma_\infty$}} \ar@/^0.6cm/[rrrr]^{\iota_x} & &
    {\widetilde{\locus(H_x)}} \ar[rr]_{\txt{\scriptsize normalization}}
    & & {\locus(H_x)} }
  $$
  
  In order to construct the divisor $\sigma_0$, choose a tangent
  vector $\vec v\in T_X|_x$ which is not contained in $F|_x$. By
  lemma~\ref{lem:special_section}, this gives a section $s(\vec v) \in
  H^0(U_x,\iota_x^*(L))$. Since $\iota_x^*(L)$ intersects the fibers
  of $\pi_x$ with mulitplicity one, and since
  $\iota_x^*(L)|_{\sigma_\infty}$ is trivial, it follows that the
  associated divisor $\Div(s(\vec v)) \in \Div(U_x)$ is a section
  $\sigma_0 \subset U_x$ which is disjoint to $\sigma_\infty$.
\end{proof}

\bibliographystyle{alpha}

\end{document}